\documentclass[epic,eepic,11pt]{amsart}

\usepackage{amsfonts,mathrsfs}
\usepackage[mathscr]{eucal}

\usepackage{amssymb}
\usepackage{amscd}
\usepackage{fancyhdr}

\usepackage[small, nohug,heads=littlevee]{diagrams}
\diagramstyle[labelstyle=\scriptstyle]

\usepackage{euler,eucal}

\DeclareMathAlphabet{\mathbf}{T1}{ppl}{bx}{n}
\DeclareMathAlphabet{\mathrm}{T1}{ppl}{m}{n}

\numberwithin{equation}{section}

\newcommand\note[1]%
{$^\dagger$\marginpar{\footnotesize{$^\dagger${#1}}}}

\def\({\left(}
\def\){\right)}
\def\<{\left<}
\def\>{\right>}

\def\newi {\sqrt{-1}\, }

\newtheorem{theorem}{Theorem}[section]
\newtheorem{proposition}[theorem]{Proposition}
\newtheorem{lemma}[theorem]{Lemma}
\newtheorem{definition}[theorem]{Definition}

\newtheorem{corollary}[theorem]{Corollary}

\theoremstyle{definition}
\newtheorem{example}[theorem]{Example}
\newtheorem{remark}[theorem]{Remark}

\newcommand\lie{\mathfrak}
\renewcommand\t{\lie{t}}
\newcommand\g{\lie{g}}

\newcommand\C{\mathbb{C}}

\renewcommand\P{\mathbb{P}}
\newcommand\CP{\mathbb{CP}}

\newcommand\J{\mathcal{J}}

\newcommand     {\comment}[1]   {}
\newcommand{\mute}[2] {}
\newcommand     {\printname}[1] {}

\newcommand\funclim[1]{\operatorname*{\mathrm{#1}}}

\renewcommand\lim{\funclim{lim}}

\newcommand\sur{\mathrel{\to\kern-1.8ex\to}}
\newcommand\iso{\mathrel{\hookrightarrow\kern-1.8ex\to}}

\newcommand\longhookrightarrow{\lhook\joinrel\longrightarrow}

\newcommand\longsur{\mathrel{\longrightarrow\kern-1.8ex\to}}
\newcommand\longiso{\mathrel{\longhookrightarrow\kern-1.8ex\to}}

\begin{document}

\bibliographystyle{amsalpha}
\date{\today}

\title{Generalized geometry, equivariant $\bar{\partial}\partial
$-lemma, and torus actions }

\author{Yi Lin}
\address{Department of Mathematics, University
of Toronto, Canada, M5S2E4, yilin@math.toronto.edu}

\begin{abstract}
In this paper we first consider the Hamiltonian action of a compact
connected Lie group on an $H$-twisted generalized complex manifold
$M$. Given such an action, we define generalized equivariant
cohomology  and generalized equivariant Dolbeault cohomology. If the
generalized complex manifold $M$ satisfies the
$\bar{\partial}\partial$-lemma, we prove that they are both
canonically isomorphic to $(S\g^*)^G\otimes H_H(M)$, where
$(S\g^*)^G$ is the space of invariant polynomials over the Lie
algebra $\g$ of $G$, and $H_H(M)$ is the $H$-twisted cohomology of
$M$. Furthermore, we establish an equivariant version of the
$\bar{\partial}\partial$-lemma, namely
$\bar{\partial}_G\partial$-lemma, which is a direct generalization
of the $d_G\delta$-lemma \cite{LS03} for Hamiltonian
 symplectic manifolds with the Hard Lefschetz property.

Second we consider the torus action on a compact generalized
K\"ahler manifold which preserves the generalized K\"ahler structure
and which is equivariantly formal.  We prove a generalization of a
result of Carrell and Lieberman \cite{CL73} in generalized K\"ahler
geometry. We then use it to compute the generalized Hodge numbers
for non-trivial examples of generalized K\"ahler structures on
 $\C\P^n$ and $\CP^n$ blown up at a fixed point.

\end{abstract}

\maketitle

\section{Introduction}

Generalized complex geometry, as introduced by Hitchin \cite{H02}
and further developed by Gualtieri \cite{Gua03}, provides a unifying
framework for both symplectic and complex geometry. It is no
surprise that many facts in complex geometry have their counterparts
in generalized complex geometry. For instance, it is well known that
a complex structure induces a $(p,q)$-decomposition for differential
forms and a splitting $d=\partial+\bar{\partial}$. Analogously,
Gualtieri \cite{Gua03} proved the presence of an $H$-twisted
generalized complex structure on a manifold determines an
alternative grading of differential forms and a similar splitting
$d_H=\partial+\bar{\partial}$, where $H$ is a closed three form and
$d_H=d-H \wedge$ is the twisted exterior derivative. (For the
twisted case, see also the appendix of \cite{KL04}.) Therefore it
makes perfect sense to define the generalized Dolbeault cohomology
and the $\bar{\partial}\partial$-lemma for $H$-twisted generalized
complex manifolds. Namely, a twisted generalized complex manifold is
said to satisfy the $\bar{\partial}\partial$-lemma if
\[ \text{ker}\partial \cap \text{im} \bar{\partial} =\text{im}\partial
\cap \text{ker} \bar{\partial} =\text{im}\bar{\partial}\partial.\]

Indeed, the $\bar{\partial}\partial$-lemma in generalized geometry
has been studied extensively by Cavalcanti in his thesis
\cite{Gil05}. It is interesting to study it for many reasons. When a
generalized complex structure is induced by a symplectic structure,
the $\bar{\partial}\partial$-lemma is equivalent to the Hard
Lefschetz property, as established by Merkulov \cite{Mer98} and
Guillemin \cite{Gu01}; whereas when a generalized complex structure
is induced by a complex structure, the
$\bar{\partial}\partial$-lemma coincides with the usual
$\bar{\partial}\partial$-lemma in complex geometry, which is known
to carry a lot of topological information. (See for instance
\cite{DGMS75}.) Recently, Gualtieri \cite{Gua04} proved that a
compact $H$-twisted generalized K\"ahler manifold satisfies the
$\bar{\partial}\partial$-lemma with respect to both generalized
complex structures involved. This result plays an important role in
the remarkable works  of \cite{li05} and \cite{Go05} which assert
that the moduli space of generalized complex structures on a compact
$H$-twisted generalized Calabi-Yau manifold is unobstructed.

In this paper we consider the consequence of the
$\bar{\partial}\partial$-lemma for group actions on generalized
complex manifolds. Sources from both symplectic and  complex
geometry have served as motivation for this work.

In \cite{LS03} Sjamaar and the author studied the Hamiltonian action
of a compact connected Lie group on a symplectic manifold with the
Hard Lefschetz property. The main results are an equivariant version
of the symplectic $d\delta$-lemma, i.e., the $d_G\delta$-lemma, and
a stronger version of Kirwan-Ginzburg formality theorem which says
that each cohomology class has a \textbf{canonical} equivariant
extension.

Motivated by \cite{LS03}, we consider the Hamiltonian action of a
compact connected Lie group on an $H$-twisted generalized complex
manifold $(M,\mathcal{J})$  as introduced in \cite{LT05}
\footnote{We note that recently there have been considerable
interests in extending quotient and reduction to the realm of
generalized complex geometry \cite{Hu05}, \cite{LT05}, \cite{SX05},
\cite{BCG05}, \cite{V05}.}. Given such an action, we introduce two
extensions of the usual equivariant Cartan complex and define their
cohomologies to be the generalized equivariant cohomology and the
generalized equivariant Dolbeault cohomology respectively. In
contrast with the usual equivariant Cartan complex, these two
extensions both contain information from moment one forms which come
up very naturally in the definition of generalized moment maps.
Assume the manifold satisfies the $\bar{\partial}\partial$-lemma, we
prove that the generalized equivariant cohomology and the
generalized equivariant Dolbeault cohomology   are both canonically
isomorphic to $(S\g^*)^G\otimes H_H(M)$, where $(S\g^*)^G$ is the
space of invariant polynomials over $\g$ and $H_H(M)$ is the
$H$-twisted cohomology of the manifold $M$. This gives an analogue
of Kirwan-Ginzburg equivariant formality theorem in generalized
geometry. In addition, we establish an equivariant version of the
$\bar{\partial}\partial$-lemma, namely, the
$\bar{\partial}_G\partial$-lemma, which is a direct generalization
of the $d_G\delta$-lemma \cite{LS03} in symplectic geometry.

We would like to mention that recently it has been found by A.
Kapustin and A. Tomasiello \cite{KT06} that the conditions used in
\cite{LT05} to define Hamiltonian actions and reductions in
generalized K\"ahler geometry are exactly the same as the physics
conditions for $(2,2)$ gauged sigma model. This has thus provided us
physics motivations to study the properties of the Hamiltonian
generalized K\"ahler manifolds as defined in \cite{LT05}. As a
consequence of the results stated in the previous paragraph, one
derives easily the first of such properties: compact Hamiltonian
twisted generalized K\"ahler manifolds always satisfy the
$\overline{\partial}_G\partial$-lemma. In view of this and the
aforementioned result of A. Kapustin and A. Tomasiello, it will be
interesting to construct more non-trivial examples of compact
Hamiltonian generalized K\"ahler manifolds. This question has been
addressed in an accompanying short note \cite{L06} which also
provides us some interesting examples for which the
$\overline{\partial}_G\partial$-lemma holds.

The second part of this paper is guided by results from K\"ahler
geometry. Historically, holomorphic vector fields on K\"ahler
manifolds have been studied by many mathematicians. Among many other
things, a famous result of Carrell and Lieberman \cite{CL73} asserts
if on a compact K\"ahler manifold $M$ there exists a holomorphic
vector field which has only isolated zeroes, then
$H^{p,q}_{\bar{\partial}}(M)=0$ unless $p=q$. Recently, assuming the
holomorphic vector field is generated by a torus action, Carrell,
Kaveh and Puppe \cite{CKP04} gave a new proof of this result. Their
method is based on equivariant Dolbeault decomposition as recently
treated by Teleman \cite{T00} and Lillywhite \cite{Lilly03}, as well
as the \textbf{localization theorem} in equivariant cohomology
theory.

We observe that the new treatment given in \cite{CKP04} could be
adapted to the case of a torus action on a compact generalized
K\"ahler manifold under certain conditions. Indeed, if we assume the
action of the torus $T$ is equivariantly formal, then a result of
Allday \cite{All04}  shows that the equivariant cohomology of the
torus action is canonically isomorphic to $S \otimes H(M)$, where
$S$ is the space of polynomials over the Lie algebra of $T$. On the
other hand, it is shown \cite{Gil05} that assuming the
$\bar{\partial}\partial$-lemma $H(M)$ will split into the direct sum
of generalized Dolbeault cohomology groups. Therefore we actually
have a version of generalized equivariant Dolbeault decomposition
for compact generalized K\"ahler manifolds. This result, together
with the localization theorem in equivariant cohomology theory,
enables us to get a generalization of the above mentioned result of
Carrell and Lieberman in generalized K\"ahler geometry.

Actually, another motivation for this piece of work is to understand
the generalized Hodge theory developed by Gualtieri \cite{Gua04} by
some concrete examples. In \cite{LT05} Tolman and the author
developed a general method which allows one to produce non-trivial
examples of generalized K\"ahler structures on many toric varieties.
We note that the classical Hodge numbers have  been long known for
toric varieties. It is thus an interesting question if one can
calculate the generalized Hodge numbers for the examples of
generalized K\"ahler structures on toric varieties discovered in
\cite{LT05}. In this article, using the analogue of Carrell and
Lieberman's result in generalized K\"ahler geometry we compute the
generalized Hodge number for non-trivial examples of generalized
K\"ahler structures on $\CP^n$ and $\CP^n$ blown up at a fixed
point.

The plan of this paper is as follows.

Section 2 goes over some basic concepts in generalized
geometry.

Section 3 presents a quick review of equivariant de Rham theory,
including a recent result of Allday \cite{All04}.

Section 4  defines the generalized equivariant cohomology and
generalized equivariant Dolbeault cohomology for Hamiltonian actions
on twisted generalized complex manifolds. Assume that the manifold
$M$ has the $\bar{\partial}\partial$-lemma, Section 4 proves that
the two cohomologies are canonically isomorphic to $(S\g^*)^G\otimes
H_H(M)$ ; moreover, it establishes an equivariant version of the
$\bar{\partial}\partial$-lemma.

Section 5 presents a generalization of Carrell and Lieberman's
result in generalized K\"ahler geometry.

Section 6 computes the generalized Hodge number for  non-trivial
examples of generalized K\"ahler structures on $\CP^n$ and $\CP^n$
blown up at a point.

\subsection*{Acknowledgement} I would like to thank Susan Tolman for
many stimulating discussions on  Hamiltonian actions and other
aspects of generalized geometry. Indeed I was led to the present
work by our joint paper \cite{LT05} which kicked the whole project
of symmetries in generalized geometry.

 I would like to thank Lisa Jeffrey and Johan Martens for their interests
 in this work and  for many helpful discussions. Finally I would like to thank
 the Mathematics Department of the University of Toronto for support and excellent
working conditions.

                   \section{ GENERALIZED COMPLEX GEOMETRY}
  Let $V$ be an $n$ dimensional vector space. There is a natural bi-linear pairing
of type $(n, n)$ which is defined by
                            \[         \langle X +\alpha , Y +  \beta \rangle = \dfrac{1}{2}(\beta (X) +
                            \alpha(Y)).\]
A \textbf{generalized complex structure} on a vector space $V$ is an
orthogonal linear map $\mathcal{J}: V\oplus V^* \rightarrow V\oplus
V^*$ such that $\mathcal{J}^2=-1$. Let $V\subset V_{\C}\oplus
V^*_{\C}$ be the $\sqrt{-1}$ eigenspace of the generalized complex
structure $\mathcal{J}$. Then $L$ is maximal isotropic and $L \cap
\overline{L} = \{0\}$. Conversely, given a maximal isotropic
$L\subset V_{\C}\oplus V^*_{\C}$  so that  $L \cap \overline{L} =
\{0\}$, there exists an unique generalized complex structure
$\mathcal{J}$ whose $\sqrt{-1}$ eigenspace is exactly $L$.

Let $\pi: V_{\C}\oplus V^*_{\C}\rightarrow V_{\C}$ be the natural
projection. The \textbf{type} of $\mathcal{J}$ is the codimension of
$\pi(L)$ in $V_{\C}$, where $L$ is the $\sqrt{-1}$ eigenspace of
$\J$.

The Clifford algebra of $V_{\C}\oplus V_{\C}^*$ acts on the space of
forms $\wedge V^*$ via
\[(X+\xi)\cdot \alpha=\iota_X\alpha+\xi\wedge\alpha.\]

Since $\mathcal{J}$ is skew adjoint with respect to the natural
pairing on $V \oplus V^* $,  $\mathcal{J}\in so(V \oplus V^*) \cong
\wedge^2(V \oplus V^*)\subset CL (V \oplus V^*)$. Therefore there is
a Clifford action of $\mathcal{J}$  on the space of  forms $\wedge
V^*$ which determines an alternative grading : $ \wedge
V_{\C}^*=\bigoplus_{\begin{subarray} k \end{subarray}}U^k,$ where
$U^k$ is the $-k\sqrt{-1}$ eigenspace of the Clifford action of
$\J$.

Let $M$ be a manifold of dimension $n$. There is a natural pairing
of type $(n,n)$ which is defined on $TM\oplus T^*M$ by
\[ \langle X+\alpha, Y_\beta \rangle
=\dfrac{1}{2}\left(\beta(Y)+\alpha(X)\right)£¬\] and which extends
naturally to $T_{\C}M\oplus T_{\C}^*M$.

For a closed three form $H$, the $H$-\textbf{twisted Courant
bracket} of $T_{\C}M\oplus T^*_{\C}M$ is defined by the identity
\[
[X+\xi,Y+\eta]=[X,Y]+L_X\eta-L_Y\xi-\dfrac{1}{2}d\left(\eta(X)-\xi(Y)\right)+\iota_Y\iota_XH.\]

The Clifford algebra of $C^{\infty}(TM\oplus T^*M)$ with the natural
pairing acts on differential forms by
\[ (X+\xi)\cdot \alpha=\iota_X\alpha+\xi\wedge \alpha.\]

A \textbf{generalized almost complex structure} on a manifold $M$ is
an orthogonal bundle map $\mathcal{J}:TM\oplus T^*M \rightarrow
TM\oplus T^*M$ such that $\mathcal{J}^2=-1$. Moreover, $\mathcal{J}$
is an $H$-\textbf{twisted generalized complex structure} if the
sections of the $\sqrt{-1}$ eigenbundle of $\mathcal{J}$ is closed
under the $H$-twisted Courant bracket. The type of $\mathcal{J}$ at
$m\in M$ is the \textbf{type} of the restricted generalized complex
structure on $T_mM$.

Let $B$ be a closed two-form on a manifold M, and consider the
orthogonal bundle map $TM\oplus T^*M\rightarrow TM \oplus T^*M$
defined by

\[ e^B=\left(\begin{matrix} 1 & 0\\ B&1 \end{matrix}\right),\]
where $B$ is regarded as a skew-symmetric map from $TM$ to $T^*M$.
This map preserves the $H$-twisted Courant bracket.  As a simple
consequence, if $\mathcal{J}$ is an $H$-twisted generalized complex
structure on M, then $\mathcal{J}' = e^B\mathcal{J} e^{-B}$ is
another $H$-twisted generalized complex structure on M, called the
B-transform of $\mathcal{J}$. Moreover, the $\sqrt{-1}$ eigenbundle
of $\mathcal{J}'$ is $e^B(L)$, so $\mathcal{J}$ and $\mathcal{J}'$
have the same type.

  Let $(M,\mathcal{J})$ be an $H$-twisted generalized complex manifold of dimension $2n$, and let
$L$ be the $\sqrt{-1}$ eigenbundle of $\mathcal{J}$. Since
$\mathcal{J}$ can be identified with a smooth section of the
Clifford bundle $CL(TM\oplus T^*M)$, there is a  Clifford action of
$\mathcal{J}$ on the space of differential forms. Let $U^k$ be the
$-k\sqrt{-1}$ eigenbundle of $\mathcal{J}$. \cite{Gua03} shows that
there is a grading of the differential forms:
   \[                       \Omega^*(M) =\Gamma( U^{-n})\oplus\cdots\oplus \Gamma(U^0) \oplus \cdots
  \oplus \Gamma(U^n) ;\]
moreover, Clifford multiplication by sections of $L$ and
$\overline{L}$ is of degree $-1$ and $1$ respectively in this
grading. This elementary fact plays an important role in Section
\ref{ddbarlemma}.

It has been shown (See e.g. \cite{Gua03} and \cite{KL04}) that the
integrability of an $H$-twisted generalized complex structure
$\mathcal{J}$ implies that
\[ d_H=d-H\wedge : \Gamma(U^k) \rightarrow \Gamma(U^{k-1})\oplus \Gamma(U^{k+1}),\]
 which gives rise to operators $\partial$ and $\bar{\partial}$
 via the projections
\[\partial: \Gamma(U^k)\rightarrow \Gamma(U^{k-1}),
 \,\,\,\bar{\partial}:\Gamma(U^k)\rightarrow \Gamma(U^{k+1}).\]
It follows that
\[\bar{\partial}^2=\partial^2=\bar{\partial}\partial+\partial\bar{\partial}=0.\]
This leads to the following definition.

\begin{definition}\label{gdolbeaut} (c.f.\cite{Gil05}) The $k$-th \textbf{generalized
Dolbeault cohomology} of $M$ is defined to be
\[ H_{\bar{\partial}}^k(M)=\text{ker}(\Gamma(U^k)
\xrightarrow{\bar{\partial}} \Gamma(U^{k+1}))\diagup
\text{im}(\Gamma(U^{k-1})\xrightarrow{\bar{\partial}}
\Gamma(U^k)).\]\end{definition}

The effect of $B$-transforms on the above grading of differential
forms and on $\bar{\partial}$, $\partial$ operators  has been
studied in \cite{Gil05}.  Though \cite{Gil05} considers only the
untwisted generalized complex structures, it is easily seen that the
same proof extends to the twisted case as well. However, note that
our sign convention for $B$-transforms differs from that of
\cite{Gil05}.

\begin{lemma}(\cite{Gil05} )\label{Btransformoperator} Let
$B$ be a closed two form and let $\mathcal{J}_B$ be the
$B$-transform of the generalized complex structure $\mathcal{J}$.
Then we have
\begin{itemize} \item [a)] $U^k_B=e^{-B}U^k$, where $U_B^k$ denotes
the $-k\sqrt{-1}$ eigenbundle of $\mathcal{J}_B$;
\item [b)] $\bar{\partial}_B=e^{-B}\bar{\partial} e^B,
\,\,\, \partial_B=e^{-B}\partial e^B$. \end{itemize}
\end{lemma}

\begin{proposition}(\cite{Gil05}) If the generalized complex manifold
$(M, \mathcal{J})$ satisfies the $\bar{\partial}\partial$-lemma,
then
\[ H_H(M)=\bigoplus_k H^k_{\bar{\partial}}(M).\]

\end{proposition}

\begin{example}  (\cite{Gua03})\label{complexgc}

\begin{itemize}\item Let $V$ be a real vector space with  a complex
structure $I$. Then the map $\mathcal{J}:V\oplus V^* \rightarrow V
\oplus V^*$ defined by \begin{equation} \label{complexstructure}
\mathcal{J}=\left(\begin{matrix} -I^* & 0
\\0 & I\end{matrix}\right) \end{equation}  is a generalized complex structure on
$V$ with the $\sqrt{-1}$ eigenspace $L=\wedge^{0,1}V_{\C} \oplus
\wedge^{1,0}V^*_{\C}$. And one easily checks that $U^k =
\bigoplus_{q-p=k}\left(\wedge^{p,q} V_{\C}^*\right)$.

\item Now let $(M,I)$ be a complex manifold. Then (\ref{complexstructure}) defines
a generalized complex structure with the $\sqrt{-1}$ eigenbundle
\[ L=\wedge^{0,1}T_{\C} M \oplus \wedge^{1,0}T^*_{\C}M.\] And the
decomposition $d:\Gamma(U^k) \rightarrow
\Gamma(U^{k+1})\oplus\Gamma( U^{k-1})$ coincides with the usual
decomposition $d = \bar{\partial}+\partial$ on the complex manifold
$(M,I)$.
\end{itemize}
\end{example}

\begin{example} (\cite{Gil05})\label{symplecticgc}
\begin{itemize}
\item Let $(V, \omega)$ be a $2n$ dimensional symplectic vector
space. Then the map $\mathcal{J}:V\oplus V^*\rightarrow V \oplus
V^*$ defined by
\begin{equation} \label{symplectic}
\mathcal{J}=\left(\begin{matrix} 0 & -\omega^{-1} \\\omega
&0\\\end{matrix} \right)\end{equation} is a generalized complex
structure on $V$. It was proved in \cite{Gil05} that
\begin{equation}\label{symplecticgrading} U^k=\{e^{i\omega}e^{\frac{\wedge}{2i}}\alpha :
\alpha \in \wedge^{n+k}V^* \}.\end{equation}

\item Now let $(M,\omega)$ be a symplectic manifold. Then (\ref{symplectic}) defines
a generalized complex structure $\J_{\omega}$ with the $\sqrt{-1}$
eigenbundle $L=\{X-\sqrt{-1}\iota_X\omega: X \in T_{\C}M\}$. It is
easy to see that (\ref{symplecticgrading}) provides a concrete
description of the alternative grading of differential forms induced
by $\J_{\omega}$. Furthermore, we have
\begin{equation}\label{symp--transform}  -2i\partial
(
e^{i\omega}e^{\frac{\wedge}{2i}}\alpha)=e^{i\omega}e^{\frac{\wedge}{2i}}(
\delta \alpha ), \,\,\,
\bar{\partial}(e^{i\omega}e^{\frac{\wedge}{2i}}\alpha)=e^{i\omega}e^{\frac{\wedge}{2i}}\left(
d\alpha \right),\end{equation} where $\delta$ is the Koszul's
boundary operator  introduced by Koszul \cite{Kos85} and studied by
Brylinski \cite{Bry88} . As a consequence, the $k$-th generalized
Dolbeault cohomology $H_{\bar{\partial}}^k(M)=H^{n-k}(M)$ as graded
vector spaces.\end{itemize}\end{example}

Let $\Omega_{\partial}(M)=\Omega(M)\cap \text{ker} \partial$. Since
$\bar{\partial}$ anti-commutes with $\partial$, $(\Omega_{\partial},
\bar{\partial})$ is a differential complex with the differential
$\bar{\partial}$. Similarly, let $H(\Omega(M),
\partial)$ be the homology of $\Omega(M)$ with respect to
$\partial$. Then $\bar{\partial}$  induces a differential on
$H(\Omega(M), \partial)$.

The following result is a simple consequence of the
$\bar{\partial}\partial$-lemma.  The proof is left as an exercise.
(c.f. [Gil05].)

\begin{proposition}\label{quasi-iso1}  Assume that the
$\bar{\partial}\partial$-lema holds. Then the $\bar{\partial}$-chain
maps in the diagram
 \[ (\Omega(M),\bar{\partial}) \leftarrow (\Omega_{\partial}(M),\bar{\partial}) \rightarrow
 H(\Omega(M),\partial) \]
 are quasi-isomorphisms, i.e., they induce isomorphisms in
 cohomology. \end{proposition}

A manifold $M$ is said to be an $H$-twisted generalized K\"ahler
manifold if it has two commuting $H$-twisted generalized complex
structures $\J_1$,$\J_2$ such that $\langle -\J_1\J_2\xi,\xi\rangle
>0$ for any $\xi \neq 0 \in C^{\infty}(T_{\C}M \oplus  T^*_{\C}M)$,
where $\langle\cdot,\cdot \rangle$ is the canonical pairing on
$T_{\C}M \oplus  T^*_{\C}M$. The following remarkable result is due
to Gualtieri.

\begin{theorem}\label{Gua'slemma}  (\cite{Gua04})  Assume that $(M,\J_1,\J_2)$ is a
  compact $H$-twisted generalized K\"ahler manifold. Then it
   satisfies the $\bar{\partial}\partial$-lemma
with respect to both $\J_1$ and $\J_2$.
\end{theorem}

\begin{example}(\cite{Gua03}) \label{Kahlerexamples}
\begin{itemize} \item  Let $(\omega,I)$ be a {\bf genuine K\"{a}hler
structure} on a manifold $M$, that is, a symplectic structure
$\omega$ and a complex structure $J$ which are {\bf compatible},
which means that $g = - \omega J$ is a Riemannian metric. By Example
\ref{symplecticgc} and \ref{complexgc}  $\omega$ and $I$ induce
generalized complex structures $\mathcal{J}_{\omega}$ and
$\mathcal{J}_I$, respectively. Moreover, it is easy to see that
$\mathcal{J}_{\omega}$ and $\mathcal{J}_I$ commute, and that
\begin{equation} \label{metric}
-\mathcal{J}_\omega \mathcal{J}_{I}=\left(
\begin{matrix} 0 & g^{-1} \\g & 0 \end{matrix} \right)\end{equation}
is a positive definite metric on $TM \oplus T^*M$. Hence
$(\mathcal{J}_\omega, \mathcal{J}_I)$ is a generalized K\"{a}hler
structure on $M$. Since $(\J_{\omega},\J_I)$ is induced by a genuine
K\"ahler structure, sometimes we will also call $(\J_{\omega},\J_I)$
a K\"ahler structure.
\item Let $(\mathcal{J}_{\omega}, \mathcal{J}_I)$ be a generalized
K\"ahler structure induced by a genuine K\"ahler structure
$(\omega,I)$, and let $B$ be a closed two form. Then
$(e^B(\J_{\omega}),e^B(\J_I))$ is also a generalized K\"ahler
structure which is said to the $B$-transform of the K\"ahler
structure $(\omega,I)$.\end{itemize} \end{example}

\section{Equivariant de Rham theory and canonical equivariant
extensions via Hodge theory}

We begin with a rapid review of equivariant de Rham theory and refer
to \cite{GS99} for a detailed account. Let $G$ be a compact
connected Lie group and let $\Omega_G(M)=(S\g^*\otimes \Omega(M))^G$
be the Cartan complex of the $G$-manifold $M$. For brevity we will
write $\Omega=\Omega(M)$ and $\Omega_G=\Omega_G(M)$. By definition
an element of $\Omega_G$ is an equivariant polynomial from $\g$ to
$\Omega$ and is called an equivariant differential form on $M$. The
bigrading of the Cartan complex is defined by
$\Omega^{ij}_G=(S^i\g^*\otimes \Omega^j)^G$. It is equipped with a
vertical differential $1 \otimes d$, which is usually abbreviated to
$d$, and the horizontal differential $d'$, which is defined by
$d'\alpha (\xi)=-\iota_{\xi}\alpha(\xi)$. Here $\iota_{\xi}$ denotes
inner product with the vector field on $M$ induced by $\xi \in \g$.
As a total complex, $\Omega_G$ has the grading
$\Omega_G^k=\bigoplus_{2i+j=k}\Omega^{ij}$ and the total
differential $d_G=d+d'$. The total cohomology $\text{ker}d_G/
\text{im}d_G$ is the de Rham equivariant cohomology $H_G(M)$. A
fundamental fact for equivariant cohomology is the following
localization theorem.

\begin{theorem}(\textbf{Localization Theorem}) Suppose a compact
connected torus $T$ acts on a compact manifold $M$. Then the kernel
of the canonical map
\[ i^*: H_T^*(M) \rightarrow H_T^*(M^T),\]
induced by the inclusion $i:M^T\rightarrow M$ is the module of
torsion elements in $H_T(M)$, where $M^T$ is the fixed point set of
the torus $T$ action. In particular, if $H_T(M)$ is a free module
over $S$, the polynomial ring over the Lie algebra $\t$ of $T$, then
$i^*$ is an injective map.

\end{theorem}

Since $\Omega^{0j}_G=(\Omega^G)^j$, the space of invariant $j$-forms
on $M$, the zeroth column of the Cartan complex $\Omega_G$ is the
invariant de Rham complex $\Omega^G$. Because $G$ is connected,
$\Omega^G$ is a deformation retract of the ordinary de Rham complex
$\Omega$. The projection $\overline{p}: \Omega_G \rightarrow
\Omega^G$, defined by $\overline{p}(\alpha)=\alpha(0)$, is a
morphism between the Cartan complex $(\Omega_G, d_G)$ and the
ordinary de Rham complex $(\Omega, d)$. The action of $G$ is called
\textbf{equivariantly formal} if $\overline{p}$ induces a surjective
map $H_G(M)\rightarrow H(M)$. \cite{All04} explained in details that
this definition of equivariant formality is equivalent to the one
that the spectral sequence of the Cartan double complex relative to
the row filtration degenerates at $E_1$ stage.

Assume the $G$ action on $M$ is equivariantly formal. Allday
\cite{All04} showed how to construct canonical equivariant
extensions using Hodge theory. Let us briefly recall his
construction here. Using a $G$-invariant Riemannian metric on $M$
one defines Hodge star operator $\ast$, adjoint operator $d^*$, the
Laplacian $\bigtriangleup$, and Green's operator $G$. (It should be
clear from the context where $G$ is the Lie group and where $G$ is
Green's operator.) Since the metric is $G$-invariant, $\ast$, $d^*$,
$\bigtriangleup$, and $G$ are all $G$-equivariant operators.
Therefore $P=(1\otimes d^*G)d'$ is a well-defined operator on the
equivariant Cartan complex. The following result is due to Allday.

\begin{theorem}(\cite{All04}) Assume the $G$-action on $M$ is
equivariantly formal. Let $\alpha \in \Omega^G$ be a closed form
(i.e.,$ d \alpha=0$). Let
\[
\widehat{\alpha}=(1-P)^{-1}\alpha=\alpha+P(\alpha)+P^2(\alpha)+\cdots+P^n(\alpha)+\cdots.
\]Then $d_G\alpha=0$.
Hence the map $\alpha \rightarrow [\widehat{\alpha}]_G$, restricted
to harmonic forms, is a canonical section of the projection
$H_G(M)\rightarrow H(M)$. \end{theorem}

As a direct consequence we have

\begin{proposition}\label{eqdecomposition}   Assume the $G$-action on $M$ is equivariantly
formal. Then there is a canonical isomorphism
\[ H_G(M)\cong (S\g^*)^G\otimes H(M).\]

\end{proposition}

\section{Generalized equivariant cohomology and the
$\bar{\partial}_G\partial$-lemma} \label{ddbarlemma}

First we recall the definition of Hamiltonian actions on $H$-twisted
generalized complex manifolds.

\begin{definition}(\cite{LT05}) \label{deftmm}
Let a compact Lie group $G$ with Lie algebra $\g$ act on a manifold
$M$, preserving an $H$-twisted generalized complex structure
$\mathcal{J}$, where $H \in \Omega^3(M)^G$ is closed. Let $L \subset
T_\C M \oplus T^*_\C M$ denote the $\sqrt{-1}$ eigenbundle of $\J$.
A {\bf twisted generalized moment map}  is a smooth function $f
\colon  M \to \g^*$ so that
\begin{itemize}
\item
There exists a one form $\eta \in \Omega^1(M,\g^*)$, called the {\bf
moment one form}, so that $\xi_M - \newi (df^\xi+  \newi \eta^\xi)$
lies in $L$ for all $\xi \in \g$, where $\xi_M$ denote the induced
vector field.
\item $f$ is equivariant.
\item $\iota_{\xi_M}H=d\eta^{\xi}$ for any $\xi \in \g$.
\end{itemize}
\end{definition}

Let $\Omega_G$ be the Cartan double complex of the $G$-manifold $M$.
Then the horizontal differential of the Cartan complex is defined by
$d'\alpha=-\iota_{\xi_M}\alpha(\xi)$, and the vertical differential
is $d$. By definition, the Cartan double complex does not encode any
information from the moment one form $\eta$ which comes up very
naturally in the definitions of generalized moment maps. Observe
that $\iota_{\xi_M}\alpha=\xi_M \cdot \alpha$, where $\cdot$ denotes
the spin action.  To extract the full information from the
Definition \ref{deftmm}, it is thus reasonable to extend $d'$ to a
new operator $\mathscr{A}$ by
\[ (\mathscr{A}\alpha)(\xi)=\mathscr{A}(\xi)\cdot
\alpha(\xi)=-\iota_{\xi_M}\alpha+\sqrt{-1}(df^{\xi}+\sqrt{-1}\eta^{\xi})\wedge
\alpha,\] where
$\mathscr{A}(\xi)=-\xi_M+\sqrt{-1}(df^{\xi}+\sqrt{-1}\eta^{\xi})$.
And now that the generalized complex manifold is $H$-twisted, it is
natural to replace the usual derivative $d$ by the twisted one
$d_H=d-H\wedge$.  Since $H$ is an invariant form, the twisted
exterior derivative $d_H$ is $G$-equivariant and induces a well
defined operator $1 \otimes d_H $ on $(S\g^* \otimes \Omega(M))^G$.
For brevity let us also denote this by $d_H$. Then we have
\[\begin{split} (d_H\mathscr{A}\alpha)(\xi)&=d_H\left(
-\iota_{\xi_M}\alpha(\xi)+\sqrt{-1}df^{\xi}\wedge \alpha(\xi)
-\eta^{\xi}\wedge \alpha(\xi)\right)\\&=
\left(-d\iota_{\xi_M}\alpha(\xi)+H\wedge
\iota_{\xi_M}\alpha(\xi)\right)-\sqrt{-1}df^{\xi}\wedge\left(d
\alpha(\xi)-H\wedge \alpha(\xi)\right) \\& \,\,\,\,\,+
\eta^{\xi}\wedge\left(d\alpha(\xi)-H\wedge
\alpha(\xi)\right)-d\eta^{\xi}\wedge \alpha(\xi)
\\&=\iota_{\xi_M}d\alpha(\xi)-L_{\xi_M}\alpha(\xi)+ H \wedge
\iota_{\xi_M}\alpha(\xi)-\sqrt{-1}df^{\xi}\wedge
d_H\alpha(\xi)+\eta^{\xi}\wedge d_H\alpha(\xi) \\&
\,\,\,\,-\iota_{\xi_M}H \wedge \alpha(\xi)\,\,\,\,(\text{ Because }
d\eta^{\xi}=\iota_{\xi_M}H.)\\&
=(\xi_M-\sqrt{-1}df^{\xi}+\eta^{\xi})\cdot
d_H\alpha(\xi)\,\,\,\,(\text{ Because } L_{\xi_M}\alpha(\xi)=0. )\\
&=\left(-\mathscr{A}d_H \alpha\right)(\xi) .\end{split} \]

This shows clearly that $d_H\mathscr{A}=-\mathscr{A}d_H$. We propose
the following definition.

\begin{definition} (\textbf{generalized equivariant cohomology})
 Let $\Omega_G=(S\g^* \otimes \Omega(M))^G$ be $Z_2$ graded. Then
$D_G=d_H+\mathscr{A}$ is a differential of degree $1$. And the $Z_2$
graded generalized equivariant cohomology is defined to be
\[ H^{\text{even/odd}}(\Omega_G,
D_G)=\dfrac{\text{ker}\left(\Omega_G^{\text{even/odd}} \xrightarrow{
D_G}
\Omega_G^{\text{odd/even}}\right)}{\text{im}\left(\Omega_G^{\text{odd/even}}
\xrightarrow{ D_G} \Omega_G^{\text{even/odd}}\right)}.\]
\end{definition}

As we are going to show in Example \ref{B-invariance}, the
generalized equivariant cohomology is invariant under $G$-invariant
$B$-transforms. This suggests that the generalized equivariant
cohomology is something natural to work with in the category of
generalized geometry. It will be interesting to define it for more
general actions and  study its property in some depth. We will leave
it for future work.

The presence of the $H$-twisted generalized complex structure
determines a splitting $d_H=\bar{\partial}+\partial$. And since
$\mathcal{J}$ is $G$ invariant, the operators $\bar{\partial}$ and
$\partial$ are $G$ equivariant. So on $(S\g^*\otimes \Omega(M))^G$
there are well-defined operators $1 \otimes \bar{\partial}$ and $1
\otimes \partial$ which we will abbreviate to $\bar{\partial}$ and
$\partial$. The following lemma says that $\bar{\partial}$ and
$\partial$ also anti-commute with the operator $\mathscr{A}$ we
introduced above.

\begin{lemma} \label{anti-commute} For any $\alpha \in (S\g^* \otimes
\Omega)^G$, we have
\[\bar{\partial}\mathscr{A}\alpha=-\mathscr{A}\bar{\partial}\alpha,\,\,\,
\partial\mathscr{A}\alpha=-\mathscr{A}\partial\alpha.\]
\end{lemma}

\begin{proof} Without the loss of generality we may assume that for
any $\xi \in \g$, $\alpha(\xi) \in U^k$, the $-k\sqrt{-1}$
eigenspace of the Clifford action of $\J$ on the space of
differential forms. Observe that
$(d_H\mathscr{A})\alpha(\xi)=-(\mathscr{A}d_H)\alpha(\xi)$, i.e.,
$d_H\mathscr{A}(\xi)\alpha(\xi)\\=-\mathscr{A}(\xi)d_H\alpha(\xi)$
for any $\xi \in \g$. Compare the $U^{k+1}$ and $U^{k-1}$ components
of $d_H\mathscr{A}(\xi)\alpha(\xi)$ and
$-\mathscr{A}(\xi)d_H\alpha(\xi)$ respectively, we conclude that for
any $\xi \in\g$,
$\mathscr{A}(\xi)\bar{\partial}\alpha(\xi)=-\bar{\partial}\alpha(\xi)$,
$\mathscr{A}(\xi)\partial\alpha(\xi)=-\partial\alpha(\xi)$.
\end{proof}

\begin{definition} (\textbf{generalized equivariant Dolbeault cohomology})
Define  $U_G=(S\g^* \otimes \Omega(M))^G$ to be the double complex
with the bigrading
\[ U_G^{i,j}=(S^i\g^* \otimes \Gamma(U^{j-i}))^G,\]
where $\Gamma(U^{j-i})$ is the $(i-j)\sqrt{-1}$ eigenspace of the
Clifford action of $\mathcal{J}$ on the space of differential forms.
It is equipped with the vertical differential $\bar{\partial}$ and
the horizontal differential $\mathscr{A}$. As a total complex, $U_G$
has the grading $U^k_G=\bigoplus_{i+j=k}U_G^{i,j}$ and the total
differential $\bar{\partial}_G=\bar{\partial}+\mathscr{A}$.The
cohomology of the total complex $(U_G, \bar{\partial}_G)$ is defined
to be the generalized equivariant Dolbeault cohomology of the
Hamiltonian action.
\end{definition}

\begin{example} \label{B-invariance}  Let $G$ act on an $H$-twisted generalized complex
manifold $(M, \mathcal{J})$ with twisted generalized moment map $f$
and moment one form $\alpha$. If $B \in \Omega^2(M)^G$, then $G$
acts on the $B$-transform of $\mathcal{J}$ with  generalized moment
map $f$ and moment one form $\alpha'$, where
$(\alpha')^{\xi}=\alpha^{\xi}+\iota_{\xi_M}B$ for all $\xi \in \g$.
\begin{itemize}
\item Let $D_G^B$ be the total differential of the generalized
equivariant double complex with respect to the Hamiltonian action on
$(M, e^B\mathcal{J}e^{-B})$. Then direct calculation shows that
$D_G^B\alpha(\xi)=D_G\alpha(\xi)+(\iota_{\xi_M}B)\wedge
\alpha(\xi)$. Since $B$ is $G$-invariant, it induces an isomorphism
\[ e^B: \Omega_G \rightarrow \Omega_G, \alpha \mapsto e^B\wedge
\alpha ;\] furthermore, we have $D_Ge^B=e^BD_G^B$. This shows
immediately that \[ H(\Omega_G, D_G)\cong H(\Omega_G, D_G^B).\]

\item Let $\bar{\partial}_G^B$ be the total differential of the
generalized equivariant Dolbeault complex with respect to the
Hamiltonian $G$-action on $(M, e^B\mathcal{J}e^{-B})$. It follows
easily from Lemma \ref{Btransformoperator} that
$\bar{\partial}_Ge^B=e^B\bar{\partial}_G^B$. As a result,
\[H(U_G, \bar{\partial}_G^B) \cong H(U_G,
\bar{\partial}_G).\]

\end{itemize}

\end{example}

Since by Lemma \ref{anti-commute} the operator $\mathscr{A}$
anti-commute with $\partial: U_G^{i,j}\rightarrow U_G^{i,j-1}$,
$\bar{\partial}_G=\bar{\partial}+\mathscr{A}$ anti-commutes with
$\partial$. So it makes sense to define the
$\bar{\partial}_G\partial$-lemma. Namely, the Hamiltonian
generalized complex manifold $M$ is said to satisfy the
$\bar{\partial}_G\partial$-lemma if and only if

\[ \text{ker}\partial \cap \text{im} \bar{\partial}_G =\text{ker}\bar{\partial}_G
\cap \text{im} \partial =\text{im}\bar{\partial}_G\partial.\]

 First let us pause for a moment to point out that when the generalized complex
 structure is induced by a symplectic structure, the $\bar{\partial}_G\partial$-lemma
is equivalent to the $d_G\delta$-lemma \cite{LS03}.

\begin{example} \label{d_Gdelta-lemma} Let $G$ act on a symplectic
manifold $(M, \omega)$ with moment map $\varPhi: M \rightarrow
\g^*$, that is, $\varPhi$ is equivariant and
$\iota_{\xi_M}\omega=d\varPhi^{\xi}$ for all $\xi \in \g$. Then $G$
also preserves the generalized complex structure
$\mathcal{J}_{\omega}$, i.e., the generalized complex structure
induced by the symplectic structure $\omega$, and $\varPhi$ is a
generalized moment map for this action with zero moment one form.

 Since the symplectic structure $\omega$ and the associated Poisson
 bi-vector $\wedge$ are $G$-invariant, the operator
 $e^{i\omega}e^{\frac{\wedge}{2i}}: \Omega(M)\rightarrow
 \Omega(M)$ extends to an operator
 $e^{i\omega}e^{\frac{\wedge}{2i}}: (S\g^*\otimes
 \Omega(M))^G\rightarrow U_G$.
Moreover, for any equivariant differential form $\alpha$ and any
$\xi \in \g$, we have:
\[\begin{split}
(\mathscr{A}e^{i\omega}e^{\frac{\wedge}{2i}}\alpha)(\xi)&=
\mathscr{A}(\xi)e^{i\omega}e^{\frac{\wedge}{2i}}\alpha(\xi)\\&
=e^{i\omega}\iota_{\xi_M} e^{\frac{\wedge}{2i}}\alpha(\xi)\\
&=e^{i\omega}e^{\frac{\wedge}{2i}}\iota_{\xi_M}\alpha(\xi)\\
&=(e^{i\omega}e^{\frac{\wedge}{2i}}d'\alpha)(\xi). \end{split}
\]
This proves that for any equivariant differential form $\alpha$,
\[\mathscr{A}e^{i\omega}e^{\frac{\wedge}{2i}}\alpha=e^{i\omega}e^{\frac{\wedge}{2i}}d'\alpha.\]

This observation, together with  (\ref{symp--transform}), shows that
\[\bar{\partial}_G
e^{i\omega}e^{\frac{\wedge}{2i}}\alpha=e^{i\omega}e^{\frac{\wedge}{2i}}d_G\alpha,\,\,\,-2i\partial
e^{i\omega}e^{\frac{\wedge}{2i}}\alpha=e^{i\omega}e^{\frac{\wedge}{2i}}\delta
\alpha ,\] where $\delta$ denotes the natural extension of the
Koszul's boundary operator to equivariant differential forms.

 Therefore the generalized equivariant cohomology group is
canonically isomorphic to the equivariant cohomology group as
$(S\g^*)^G$-modules. Furthermore, it is easy to see that the
$\bar{\partial}_G\partial$-lemma is equivalent to the
$d_G\delta$-lemma \cite{LS03} which asserts that
\[ \text{ker}d_G \cap \text{im} \delta =\text{ker}\delta
\cap \text{im} d_G =\text{im} d_G\delta.\]
\end{example}

 Now that $\bar{\partial}$ anti-commutes with $\partial$, it is
 straightforward to check that $U_{G,\partial}=U_G\cap
 \text{ker}(\partial)$ is a sub-double complex of $U_G$ and that the
 homology complex $H(U_G, \partial)$ of $U_G$ with respect to
 $\partial$ is a double complex with differentials induced by
 $\bar{\partial}$ and $\mathscr{A}$. Thus we have the following
 diagram of morphisms of double complexes

 \begin{equation}  \label{morphisms} (U_G, \bar{\partial}_G)\leftarrow (U_G\cap
 \text{ker}(\partial),\bar{\partial}_G) \rightarrow H(U_G,
 \partial), \end{equation}

Since $\partial$ does not act on the polynomial part, these
morphisms are linear over the invariant polynomials $(S\g^*)^G$. Let
us first examine the homology complex $H(U_G, \partial)$. We will
need two preliminary results.

\begin{lemma} \label{trivialaction1} Suppose the action of the connected compact Lie group
$G$ on $M$ preserves the $H$-twisted generalized complex structure
$\mathcal{J}$, and suppose that $(M, \mathcal{J})$ satisfies the
$\bar{\partial}\partial$-lemma, then the induced $G$ action on
$H_{\bar{\partial}}^*(M)$ and $H^*(\Omega(M),\partial)$ is
trivial.\end{lemma}

\begin{proof} Let $\alpha$ be a representative of an element
$[\alpha]$ of $H_{\bar{\partial}}(M)$. By Proposition
\ref{quasi-iso1} we may well assume that $\alpha$ is both $\partial$
and
 $\bar{\partial}$ closed. In particular, this implies that
 $\alpha$ is $d$-closed. Let $g$ be an element of $G$. Since the
 induced Lie group action on the de Rham  cohomology is trivial, we
 have \begin{equation} \label{trivialaction2} g^*\alpha-\alpha=d\gamma \end{equation}
 for some $\gamma \in \Omega(M)$. Without the loss of generality we
 may assume that $\gamma$ has only $U^{k-1}$ component
 $\gamma^{k-1}$ and $U^{k+1}$ component $\gamma^{k+1}$. By comparing the
components of the both sides of  (\ref{trivialaction2}) we get
\[g^*\alpha-\alpha=\bar{\partial}\gamma^{k-1}+\partial\gamma^{k+1},
\,\,\partial\gamma^{k-1}=0,\,\,\bar{\partial}\gamma^{k+1}=0.\] Note
that $\partial\gamma^{k+1}$ is both $\partial$ exact and
$\bar{\partial}$ closed. By the $\bar{\partial}\partial$-lemma
$\partial\gamma^{k+1}=\bar{\partial}\partial\eta$ for some $\eta \in
U^k$. Therefore
$g^*\alpha-\alpha=\bar{\partial}(\gamma^{k-1}+\partial\eta)$. This
shows that the induced $G$ action on $H_{\bar{\partial}}(M)$ is
trivial. A similar argument shows that the induced $G$ action on
$H(\Omega(M), \partial)$ is trivial.

\end{proof}

It is important to notice that $\partial$ is not a derivation. But
we have the following Leibniz rule.

\begin{lemma}\label{leibniz} Let $f$ be the generalized moment map and let $\eta \in \Omega^1(M, \g^*)$ be the associated
moment one form. For any $\xi \in \g $, define
$\mathscr{A}(\xi)=-\xi+\sqrt{-1}(df^{\xi}+\sqrt{-1}\eta^{\xi})$ as
before. Then for any $\alpha \in \Omega(M)$, we have
\[\partial(f^{\xi}\alpha)=-\dfrac{\sqrt{-1}}{2}\mathscr{A}(\xi)\cdot
\alpha+f^{\xi}\partial\alpha.\]

\end{lemma}

\begin{proof}  Without the loss of generality we may assume that
$\alpha \in U^k$. First we note that
\[d(f^{\xi}\alpha)=df^{\xi}\wedge\alpha+f^{\xi}d\alpha.\]
It is easily seen that
$df^{\xi}=-\dfrac{\sqrt{-1}}{2}\mathscr{A}(\xi)+\dfrac{\sqrt{-1}}{2}\overline{\mathscr{A}(\xi)}$
with $-\dfrac{\sqrt{-1}}{2}\mathscr{A}(\xi) \in C^{\infty}(L)$ and
$\dfrac{\sqrt{-1}}{2}\overline{\mathscr{A}(\xi)}\in
C^{\infty}(\overline{L})$. Thus
\[d(f^{\xi}\alpha)=\left(-\dfrac{\sqrt{-1}}{2}\mathscr{A}(\xi)+\dfrac{\sqrt{-1}}{2}\overline{\mathscr{A}(\xi)}\right)\cdot \alpha+f^{\xi}d\alpha.\]
Now compare the $U^{k-1}$ component of the both side of the above
equality, we get
that\[\partial(f^{\xi}\alpha)=-\dfrac{\sqrt{-1}}{2}\mathscr{A}(\xi)\cdot
\alpha+f^{\xi}\partial\alpha.\]

\end{proof}

We note that as an immediate consequence of Lemma
\ref{trivialaction1} and Lemma \ref{leibniz}, the homology complex
$H(U_G,\partial)$ is a double complex with trivial differentials if
the manifold $M$ satisfies the $\bar{\partial}\partial$-lemma.

\begin{lemma} \label{trivialcomplex}If $(M,\mathcal{J})$ satisfies the
$\bar{\partial}\partial$-lemma, then both differentials
$\bar{\partial}$ and $\mathscr{A}$  on $H(U_G,
\partial)$ are zero.
\end{lemma}

\begin{proof}  First we observe that the (ordinary) $\bar{\partial}\partial$-lemma
holds for equivariant forms as well as for ordinary forms. The
reason is that $\partial$ and $\bar{\partial}$ acts on $U_G$ as $1
\otimes \partial$ and $1 \otimes \bar{\partial}$ respectively and
the both operators are $G$ equivariant. Now suppose $\alpha \in U_G$
satisfies that $\partial \alpha=0$. Then
$\bar{\partial}\alpha=\bar{\partial}\partial\beta
=-\partial\bar{\partial}\beta$ for some $\beta \in U_G$. Hence the
differential on $H(U_G, \partial)$ induced by $\bar{\partial}$ is
zero.

To prove the other differential is zero we have to be more careful
to pick a representative of an element of $H(U_G, \partial)$. By
Lemma \ref{trivialaction1} the induced $G$ action on $H(\Omega(M),
\partial)$ is trivial. This implies that
\begin{equation} \label{homologycomplex}
H(U_G,\partial)=\left(S\g^* \otimes H(\Omega(M),\partial) \right)^G=
(S\g^*)^G\otimes H(\Omega(M), \partial).\end{equation}

Choose a basis $\xi_i$ of $\g$. Let $x_i$ be the dual basis of
$\g^*$ and let $f_i$ be a basis of the vector space $(S\g^*)^G$ of
invariant polynomials. It follows from  (\ref{homologycomplex}) that
an element of $H(U_G, \partial)$ can be represented by an $\alpha
\in U_G$ with $\partial\alpha=0$ of the form
$\alpha=\sum_if_i\otimes \alpha_i$ for unique $\alpha_i\in
\Omega^G(M)$. It follows that $\partial\alpha_i=0$ for all $i$. And
so by Lemma \ref{leibniz}
$\mathscr{A}(\xi_j)\alpha_i=\partial\beta_{ij}$. Hence
$\mathscr{A}\alpha=\sum_{i,j}x_jf_i\otimes
\partial\beta_{ij}=\partial\left( \sum_{i,j}x_jf_i\otimes
\beta_{ij}\right)$. Since the operator $\mathscr{A}$ and $\partial$
are equivariant, after averaging over $G$ we get
$\mathscr{A}\alpha=\partial\beta$ with $\beta \in U_G$, i.e., the
differential on $H(U_G, \partial)$ induced by $\mathscr{A}$ is
trivial.

\end{proof}

Let $E$ be the spectral sequence of $U_G$ relative to the filtration
associated to the horizontal grading and $E_{\partial}$ that of
$U_{G,\partial}$. The first terms are
\begin{equation} \label{E1terms} \begin{split}
& E_1=\text{ker}\bar{\partial}/ \text{im}\bar{\partial}=\left(S\g^*
\otimes H_{\bar{\partial}}(M)\right)^G=(S\g^*)^G\otimes
H_{\bar{\partial}}(M),
\\&(E_{\partial})_1=(\text{ker}\bar{\partial}\cap
\text{ker}\partial)/(\text{im}\bar{\partial}\cap
\text{ker}\partial)=(S\g^* \otimes
H(\Omega_{\partial}(M),\bar{\partial}))^G=(S\g^*)^G\otimes
H_{\bar{\partial}}(M).
\end{split}\end{equation}

Here we used the isomorphism
$H(\Omega_{\partial}(M),\bar{\partial})\cong H_{\bar{\partial}}(M)$
of Proposition \ref{quasi-iso1} and the connectedness of $G$. By
Lemma \ref{trivialcomplex} $H(U_G,
\partial)$ is a double complex with trivial differentials, so its spectral sequence is
constant with trivial differential at each stage. The two morphisms
(\ref{morphisms}) induce morphisms of spectral sequences
\[ E \leftarrow E_{\partial} \rightarrow H(U_G, \partial).\]

It follows from (\ref{homologycomplex}) and (\ref{E1terms}) that
these two morphisms are isomorphisms at the first stage and hence
are isomorphisms at every stage. So they induce isomorphisms on the
total cohomology. In fact, since the spectral sequence for $H(U_G,
\partial)$ is constant, so are spectral sequences $E$ and
$E_{\partial}$. This proves the following result, where $H_{G,
\partial}(M)$ denotes the total cohomology of $U_{G,\partial}$.

\begin{theorem} (\textbf{equivariant formality I})
\label{eqformality1} Assume that the generalized complex manifold
$M$ satisfies the $\bar{\partial}\partial$-lemma. Then the spectral
sequences $E$ and $E_{\partial}$ degenerate at the first terms. And
the morphisms (\ref{morphisms}) induce isomorphisms of
$(S\g^*)^G$-modules
\[ H(U_G, \bar{\partial}_G)\leftarrow H_{G,
\partial}(M)\rightarrow (S\g^*)^G\otimes H_{\partial}(M).\]

\end{theorem}

The following corollary is an immediate consequence of  Theorem
\ref{eqformality1} and Proposition \ref{quasi-iso1}.

\begin{corollary} \label{eqformality2} Assume that the generalized complex manifold $M$
satisfies the $\bar{\partial}\partial$-lemma. As $(S\g^*)^G$-modules
\[ H(U_G, \bar{\partial}_G)\cong (S\g^*)^G\otimes H_H(M).\]

\end{corollary}

To prove the $\bar{\partial}_G\partial$-lemma we need the following
useful technical lemma. For a proof, we refer to \cite{L04}.

\begin{lemma}(\cite{L04}) \label{Ddeltalemma} ($D\delta$-lemma) Let $(K^{**}, d, d')$ be
a double complex which is bounded in the following sense: for each
$n$, there are only finitely many non-zero components in the direct
sum $K^n=\bigoplus_{i+j=n}K^{i,j}$. Here $d$ is the degree $1$
vertical differential and $d'$ the degree $1$ horizontal
differential. Assume that there is a degree $-1$ vertical
differential $\delta$ which anti-commutes with both $d$ and
$d'$,i.e., $d\delta=-\delta d$, $d'\delta=-\delta d'$. Let $(K^*,
D)$ be the associated total complex, where $D=d+d'$. And assume that
the double complex $(K^{**}, d,d')$ satisfies:
\begin{itemize}\item [a)]
$\text{im}d\cap\text{ker}\delta=\text{ker}d \cap
\text{im}\delta=\text{im}d\delta$;

\item [b)] The spectral sequence associated to the row filtration
degenerates at the $E_1$ stage. \end{itemize} Then we have
$\text{im}D \cap \text{ker}\delta=\text{im}D\delta$.
\end{lemma}

We are ready to state the  main result of this section.

\begin{theorem}\label{eqddbar-lemma} Consider the Hamiltonian action
of a connected compact Lie group on an $H$-twisted generalized
complex manifold $(M,\mathcal{J})$. If $(M,\J)$ satisfies the
$\bar{\partial}\partial$-lemma, then
\[ \text{im}\bar{\partial}_G \cap
\text{ker}\partial =\text{ker}\bar{\partial}_G\cap
\text{im}\partial=\text{im} \bar{\partial}_G\partial.\]

\end{theorem}

\begin{proof} As a direct consequence of Theorem \ref{eqformality1}
and Lemma \ref{Ddeltalemma}, we have $\text{im}\bar{\partial}_G \cap
\text{ker}\partial =\text{im} \bar{\partial}_G\partial$. The second
half of the $\bar{\partial}_G\partial$-lemma follows from the first:
assume that $\bar{\partial}_G\alpha=0$ and $\alpha$ is $\partial$
exact. Then the cohomology class of $\alpha$ in $H(U_G,
\partial)$ is zero, so by Theorem \ref{eqformality1} the cohomology
class of $\alpha$ in $H_{G, \partial}(M)$ is zero, i.e., $\alpha$ is
$\bar{\partial}_G$ exact. Hence
$\alpha=\bar{\partial}_G\partial\beta$ for some $\beta$.
\end{proof}

As an easy consequence, the $\overline{\partial}_G\partial$-lemma
holds for compact $H$-generalized K\"ahler manifolds. To state this
result more precisely, let us first recall the definition of
Hamiltonian actions on twisted generalized K\"ahler manifolds.

\begin{definition} (\cite{LT05}) Let the compact Lie group $G$ with Lie algebra $\g$ act on a
manifold $M$. A {\bf generalized moment map}  for an invariant
$H$-twisted generalized K\"{a}hler structure $(\J_1,\J_2)$ is a
generalized moment map for the generalized complex structure $\J_1$.
If such a generalized moment map exists, the action on the
$H$-twisted generalized K\"ahler manifold $(M,\J_1,\J_2)$ is said to
be Hamiltonian.\end{definition}

\begin{corollary} Assume that the action of the compact Lie group $G$  on
an $H$-twisted generalized K\"ahler manifold $(M, \J_1,\J_2)$ is
Hamiltonian. Then the $\overline{\partial}_G\partial$-lemma holds
for the generalized complex manifold $(M, \J_1)$.\end{corollary}

\begin{proof} By Theorem \ref{Gua'slemma}, the
$\overline{\partial}\partial$-lemma holds on $M$ with respect to
both $\J_1$ and $\J_2$. By definition, the action of $G$ is
Hamiltonian on the generalized complex manifold $(M,\J_1)$. Now the
corollary follows easily from Theorem \ref{eqddbar-lemma}.

\end{proof}

To conclude this section let us present an application of the
$\bar{\partial}_G\partial$-lemma which says that the generalized
equivariant cohomology is canonically isomorphic to
$(S\g^*)^G\otimes H_H(M)$ provided the manifold $M$ satisfies the
$\bar{\partial}\partial$-lemma.

Observe that the inclusion map
\[ (U_{G, \partial}, \bar{\partial}_G) \hookrightarrow
(\Omega_G, D_G)\] is actually a chain map with respect to the
differentials $\bar{\partial}_G$ and $D_G$ since
$\bar{\partial}_G\alpha=D_G\alpha$ for any $\alpha \in U_{G,
\partial}=U_G\cap \text{ker}\partial$. So it induces a map
\begin{equation} \label{morphism2} H_{G, \partial}(M)\rightarrow
H(\Omega_G, D_G).\end{equation}

Suppose $\alpha$ is a representative of a cohomology class in
$H_{G,\partial}(M)$ and $\alpha=D_G\beta$ for some $\beta \in
\Omega_G$. Then $\alpha-\partial \beta=\bar{\partial}_G\beta$ is
both $\partial$ closed and $\bar{\partial}_G$ exact. So by the
$\bar{\partial}_G\partial$-lemma
 $\alpha-\partial \beta=\bar{\partial}_G\partial \gamma$ for
some $\gamma \in U_G$. Thus $\alpha=\partial(\beta-\bar{\partial}_G
\gamma)$ is both $\partial$-exact and $\bar{\partial}_G$ closed.
Applying the $\bar{\partial}_G\partial$-lemma again, we conclude
that $\alpha=\bar{\partial}_G\partial \eta$ for some $\eta \in U_G$.
This shows that $\alpha$ represents a trivial cohomology class in
$H_{G, \partial}(M)$ and the map (\ref{morphism2}) is injective. Now
suppose $\alpha$ is a representative of a cohomology class in
$H(\Omega_G, D_G)$. Then since $\bar{\partial}_G\partial
\alpha=-\partial D_G\alpha=0$, $\partial\alpha$ is both $\partial$
exact and $\bar{\partial}_G$ closed.So
$\partial\alpha=\bar{\partial}_G\partial\beta$ for some $\beta \in
U_G$. It follows that
$\partial(\alpha+\bar{\partial}_G\beta)=\partial(\alpha+D_G\beta)=0$.
Since $D_G\alpha= \bar{\partial}_G\alpha+\partial \alpha=0$,
$\bar{\partial}_G(\alpha+D_G\beta)=-\partial \alpha+
\bar{\partial}_G\partial\beta=0$. This shows clearly that the
cohomology class of $\alpha$ in $H(\Omega_G, D_G)$ is the image of
the cohomology class of $\alpha+D_G\beta$ in $H_{G,\partial}(M)$.
Hence the map (\ref{morphism2}) is surjective. The above discussion,
together with Corollary \ref{eqformality2}, leads to the following
theorem.

\begin{theorem} (\textbf{equivariant formality II}) Assume the generalized complex manifold $M$ satisfies the
$\bar{\partial}\partial$-lemma. Then
\[H(\Omega_G, D_G)\cong (S\g^*)^G\otimes H_H(M).\]

\end{theorem}

\section{Torus actions on generalized K\"ahler manifolds}

Assume that $(M, \mathcal{J})$ is a generalized complex manifold
which satisfies the $\bar{\partial}\partial$-lemma. By Theorem
\ref{quasi-iso1} we have that
$H(M)=\bigoplus_kH^k_{\bar{\partial}}(M)$. Therefore we have the
following decomposition of $H_G(M)$.

\begin{proposition}\label{eqdecomposition} Assume  the connected compact Lie group $G$
action on $M$ is equivariantly formal and assume $M$ satisfies the
$\bar{\partial}\partial$-lemma. Then there is a canonical
isomorphism of $(S\g^*)^G$-modules
\begin{equation}\label{eqdom} H_G(M)=\bigoplus_k (S\g^*)^G\otimes
H^k_{\bar{\partial}}(M).\end{equation} \end{proposition}

\begin{remark}It is obvious that the above canonical isomorphism depends only on
the invariant metric that we choose. When the generalized complex
structure $\J$ is induced by a complex structure $I$ in an
$G$-invariant K\"ahler pair $(\omega,I)$, it is easy to recover from
(\ref{eqdom}) the equivariant Dolbeault decomposition as treated by
Teleman \cite{T00} and Lillywhite \cite{Lilly03}. Indeed, assuming
the group action is equivariantly formal, using the
$\bar{\partial}\partial$-lemma for compact K\"ahler manifolds it is
not difficult to show directly that, for any $\bar{\partial}$-closed
differential form $\alpha$,  any holomorphic vector $Z$ generated by
the group action and any $p>0$, there exist differential forms
$\alpha_1,\cdots,\alpha_p$ so that
\[ \iota_Z\alpha=\bar{\partial}\alpha_1,
\cdots,\iota_Z\alpha_{p-1}=\bar{\partial}\alpha_p.\] Then one can
apply Allday's argument in \cite{All04} to the invariant K\"ahler
metric and prove that the right hand side of (\ref{eqdom}) is
canonically isomorphic to the equivariant Dolbeault cohomology. This
approach will give us a new Hodge theoretic proof of the usual
equivariant Dolbeault decomposition without using any Morse theory.
\end{remark}

 Henceforth we will assume that $G=T$ is a connected compact torus and that
 the $T$ action preserves the
generalized complex structure $\mathcal{J}$. Our next observation is
that the fixed point submanifold of the $T$ action is a generalized
complex submanifold in the sense specified in \cite{BB03}.  However,
for the convenience of the reader, let us first review the notion of
generalized complex submanifolds \cite{BB03}.

Let $W$ be a submanifold of the generalized complex manifold $(M,
\mathcal{J})$ and let $L$ be the $\sqrt{-1}$ eigenbundle of
$\mathcal{J}$. Then at each point $x\in N$ define \[
L_{W,x}=\{X+(\xi \mid_{T_{\C}W}) : X+\xi \in
L\cap\left(T_{\C,x}W\oplus T_{\C,x}^*M\right)\}.\]

This actually defines a Dirac structure on $W$, i.e., a maximal
isotropic distribution $L_W \subset T_{\C}W\oplus T_{\C}^*W$  whose
sections are closed under the Courant bracket. If $L_W$ is such that
$L_W\cap \overline{L_W}=0$, then $W$ is said to be a
\textbf{generalized complex submanifold} of $M$. It is clear from
the definition that if $W$ is a generalized complex submanifold of
$M$, then there is a unique generalized complex structure
$\mathcal{J}_W$ on $W$ whose $\sqrt{-1}$ eigenbundle is exactly
$L_W$. Moreover, \cite{BB03} gives a simple condition, the "split
condition", to ensure the submanifold $W$ of $M$ is a generalized
complex submanifold. Specifically, $W$ is said to be \textbf{split}
if there exists a smooth subbundle $N$ of $TM\mid_W$ such that
$TM\mid_W=TW\oplus N$ and such that $TW\oplus \text{Ann}( N)$ is
invariant under the generalized complex structure $\mathcal{J}$.
Here $\text{Ann}(N) \subset T^*M\mid_W$ denotes the annihilator of
$N$.

\begin{proposition}\label{split}(\cite{BB03}) Let $(M, \mathcal{J})$
be a generalized complex manifold and let $W$ be a split
submanifold. Then $W$ is a generalized complex submanifold of $M$.
Moreover, let $\psi: TW\oplus\text{Ann}(N) \rightarrow TW\oplus
T^*W$ be the natural isomorphism, then the induced generalized
complex structure $\mathcal{J}_W$ on $W$ has the form
$\mathcal{J}_W=\psi\circ \mathcal{J} \circ \psi^{-1}$.

\end{proposition}

Recall that the fixed point submanifold of a symplectic torus action
on a symplectic manifold is a symplectic submanifold. The following
lemma is a generalization of this well-known fact. We note that in
the case of  $Z_2$ actions on generalized complex manifolds the
similar result has been obtained by J. Barton and M.  Sti\'{e}non
\cite{BS06} using different methods.

\begin{lemma}\label{fixedpointsplit} Suppose that the torus $T$ acts on
the generalized complex manifold $(M, \mathcal{J})$ preserving the
generalized complex structure $\mathcal{J}$. And suppose that $Z$ is
a connected component of the fixed point submanifold. Then $Z$ is a
split submanifold and so is a generalized complex submanifold of $M$

\end{lemma}

\begin{proof} Let $x \in Z$ be a fixed point of the torus action.
Then the action of the torus $T$ on $M$ induces a $T$-module
structure on $T_xM$ and a dual $T$-module structure on $T_x^*M$. Let
$\{\vartheta_1,\vartheta_2,\dots,\vartheta_k\}$ be the set of all
distinct weights of the $T$-module $T_xM$, where $\vartheta_1 \equiv
1$. Then
$\{\vartheta^{-1}_1,\vartheta^{-1}_2,\dots,\vartheta^{-1}_k\}$ is
the set of all  distinct weights of the dual $T$-module $T^*_xM$.
Let $V_i \subset T_xM$ be the weight space corresponding to
$\vartheta_i$ and $V^*_i\subset T_x^*M$ the weight space
corresponding to $\vartheta_i^{-1}$, $1 \leq i \leq k$, and let
$N_x=\bigoplus_{i=2}^kV_i$. Then it is clear  $V_1=T_xZ$. Moreover,
we claim that $\text{Ann}(N_x)$, the annihilator of $N_x$, coincides
with $V_1^*$. Indeed, for any $u^*\in V_1^*$, $w\in V_i$,$i\geq 2$,
we have
\[ \langle w, u^*\rangle =\langle tw, t u^*\rangle= \langle \vartheta_i(t)w,u^*
\rangle=\vartheta_i(t) \langle w,u^*\rangle,\] where $t$ is an
arbitrarily chosen element in $ T$. It follows  $\langle u^*,
w\rangle=0$ and so $V_1^* \subset \text{Ann}(N_x)$. A dimension
count shows that we actually have $V_1^*=\text{Ann}(N_x)$. Since the
torus action preserves the generalized complex structure $\J$, $V_1
\oplus V^*_1=T_xZ \oplus \text{Ann}(N_x)$ is invariant under $\J$.
Put $N=\bigcup\limits_{x}   N_x$. Then $N$ is a smooth subbundle of
$TM\mid_Z$ such that $TM\mid_Z=TZ\oplus N$ and such that $TZ\oplus
\text{Ann}(N)$ is invariant under $\J$. This completes the proof
that $Z$ is a split submanifold.
\end{proof}

\begin{proposition}\label{fixedpoint}  Suppose that the torus $T$ action preserves the
generalized K\"ahler structure $(\J_1,\J_2)$. And suppose $Z$ is a
connected component of the fixed point submanifold of the torus
action. Then $Z$ is split with respect to both $\J_1$ and $\J_2$.
Furthermore, the induced generalized complex structures $\J_{1,Z}$
and $\J_{2,Z}$ defines a generalized K\"ahler pair on $Z$.

\end{proposition}
\begin{proof}The first assertion is an immediate consequence of Lemma
\ref{fixedpointsplit}. The second assertion follows from the
description of the induced generalized complex structures on a split
submanifold given in Proposition \ref{split}.
\end{proof}

\begin{example} As explained in \cite{BB03}, if $\J_{\omega}$ is a
generalized complex structure induced by a symplectic structure
$\omega$ on $M$, then $Z$ is a generalized complex submanifold of
$(M,\J_{\omega})$ if and only if $Z$ is a symplectic submanifold of
$(M, \omega)$. By Example \ref{symplecticgc}
\[\Gamma\left(U^k(M)\right)=\{e^{i\omega}e^{\frac{\wedge}{2i}}\alpha :\alpha\in
\Omega^{n+k}(M)\},\,\,\,\Gamma\left(U^k(Z)\right)=\{e^{i\omega_0}e^{\frac{\wedge_0}{2i}}\alpha
:\alpha\in \Omega^{n+k}(Z)\},\]where $\omega_0$ is the restriction
of $\omega$ to $Z$ and $\wedge_0$ is the Poisson bivector on $Z$
associated to the symplectic structure $\omega_0$. It is then not
hard to see that $\alpha \in \Gamma\left(U^k(M)\right)$ does not
necessarily imply that $(\alpha\mid_Z)
\in\Gamma\left(U^k(Z)\right)$. As a result, we see
$\bar{\partial}\alpha=0$ on $M$ does not necessarily imply that
$\bar{\partial}_Z(\alpha\mid_Z)=0$, where $\bar{\partial}_Z$ is
$\bar{\partial}$ operator associated to the generalized complex
structure on $Z$ induced by the symplectic structure $\omega_0$ .

\end{example}

\begin{example} \label{complexsubmanifold} \begin{itemize}
 \item Let $V$ be a real vector space
with a complex structure $I$, and let $W$ be a subspace of $V$ which
is invariant under $I$. It is clear that $W$ inherits a complex
structure $I_W$ from $(V,I)$. Denote by $\J$ and $\J_W$ the
generalized complex induced by $I$ and $I_W$ respectively. Let
$U^k(V)$ and $U^k(W)$ be the $-k\sqrt{-1}$ eigenspace of the
Clifford actions of $\J$ and $\J_W$ respectively. Then by Example
\ref{complexgc} \[ U^k(V)=\bigoplus_{p-q=k}\wedge^{p,q}V^*,\,\,\,
U^k(W)=\bigoplus_{p-q=k}\wedge^{p,q}W^*.\] In particular, if $\alpha
\in U^k(V)$, then $(\alpha \mid_W ) \in U^k(W)$.
\item  Suppose that $\J$ is a generalized complex structure on $M$
induced by a complex structure $I$. Then as shown in \cite{BB03} a
submanifold $Z$ is a generalized complex submanifold of $(M,\J)$ if
and only if $S$ is a complex submanifold of $(M,I)$.

Now suppose that $Z$ is a generalized complex submanifold of
$(M,\J)$ with the induced generalized complex structure $\J_W$. By
the preceding discussion, if $\alpha \in \Gamma\left(U^k(M)\right)$,
then $(\alpha \mid_Z) \in \Gamma\left(U^k(Z)\right)$. So if $\alpha
\in \Gamma\left(U^k(M)\right)$ such that $\bar{\partial}\alpha=0$,
then $\bar{\partial}_Z(\alpha \mid_Z)=0$, where $\bar{\partial}_Z$
is the $\bar{\partial}$-operator associated to the generalized
complex structure $\J_W$.
\end{itemize} \end{example}

The same argument  actually gives us the following slightly more
general result.

\begin{lemma}\label{restriction} Let $(M,\J)$ be a generalized complex manifold and
let $Z$ be a generalized complex submanifold with the generalized
complex structure $\J_Z$. Suppose that for each $x \in Z$, $\J_x$ is
induced by a complex structure on the tangent space $T_xM$. As a
result, the generalized complex structure $\J_Z$ on $Z$ is induced
by a complex structure on $Z$. Furthermore if $\alpha \in
\Gamma\left(U^k(M)\right)$, then $(\alpha \mid_Z) \in
\Gamma\left(U^k(Z)\right)$. In particular, this implies that if
$\bar{\partial}\alpha=0$, then $\bar{\partial}_Z(\alpha \mid_Z)=0$,
where $\bar{\partial}_Z$ is the $\bar{\partial}$-operator associated
to the generalized complex structure $\J_W$.
\end{lemma}

\begin{theorem} \label{gcarrell}Consider the action of a torus $T$ on a generalized
K\"ahler manifold $M$ which preserves the generalized K\"ahler
structure $(\J_1,\J_2)$. Assume that the action is equivariantly
formal. Let $Z$ be the fixed point submanifold. And assume that  for
any $x \in Z$, $\J_{2,x}$, the restriction of $\J_2$ to $T_xM$,   is
induced by a complex structure on $T_xM$. Then
\[ H^i_{\bar{\partial}_2}(M)=0\,\,\,\text{if }\vert i\vert
>\text{dim}Z,\]
where $\bar{\partial}_2$ is the $\bar{\partial}$-operator associated
to the generalized complex structure $\J_2$.

\end{theorem}

\begin{proof} Let us denote by $S$ the polynomial ring over the Lie
algebra $\t$ of $T$. By Proposition \ref{eqdecomposition} there is a
canonical isomorphism $H_T(M)\cong \bigoplus_{k}S\otimes
H_{\bar{\partial}}^k(M)$. By Proposition \ref{fixedpoint} the fixed
point submanifold $Z$ is a compact generalized K\"ahler manifold and
therefore satisfies the $\bar{\partial}\partial$-lemma with respect
to both induced generalized complex structures. This implies that
$H_T(Z)\cong\bigoplus_kS\otimes H^k_{\bar{\partial}}(M)$. A
straightforward check shows that there is a commutative diagram
\[ \begin{CD}H_T(M) @>{\operatorname{\text{isomorphism}}}>> \bigoplus_{k}S\otimes
H_{\bar{\partial}}^k(M)\\
@V{\operatorname{injection}}VV  @VVV\\
H_{T}(Z)@>{\operatorname{\text{isomorphism}}}>>
\bigoplus_{k}S\otimes H_{\bar{\partial}}^k(Z),
\end{CD}\]  where the right vertical map
is well-defined because of Lemma \ref{restriction}. Observe that the
direct summand $S\otimes H^k_{\bar{\partial}}(M)$ is mapped into
$S\otimes H^k_{\bar{\partial}}(Z)$; moreover, since the above
diagram is commutative, the map $S\otimes
H^k_{\bar{\partial}}(M)\rightarrow S\otimes H^k_{\bar{\partial}}(Z)
$  is injective. By Lemma \ref{restriction} the generalized complex
structure on $Z$ is induced by a complex structure and so
$\Gamma\left(U^k(Z)\right)=0$ if $\vert i \vert
>\text{dim}Z$ by type consideration. Therefore for any $\vert i \vert
>\text{dim}Z$,
$H^k_{\bar{\partial}}(Z)=0$ and so $H^k_{\bar{\partial}}(M)=0$.

\end{proof}

When the invariant generalized K\"ahler structure $(\J_1,\J_2)$ is
induced by an invariant K\"ahler structure $(\omega,I)$, it is easy
to recover from Theorem \ref{gcarrell} the following result of
Carrell and Lieberman \cite{CL73}, \cite{CKP04}.

\begin{corollary} Suppose that $M$ is a compact K\"ahler manifold with an
equivariantly formal torus  action which preserves the K\"ahler
structure. And suppose that the fixed point submanifold $Z$ of the
torus action  is non-empty.  Assume that $\vert p-q \vert >
\text{dim} Z$. Then $ H_{\bar{\partial}}^{p,q}(M)=0.$
\end{corollary}

\begin{proof}  Since the fixed point of the torus action is
non-empty, by a well-known result of Frankel \cite{F59} the torus
action is Hamiltonian and so is equivariantly formal by
Kirwan-Ginzburg equivariant formality theorem. Let $(\J_1,\J_2)$ be
the generalized K\"ahler structure induced by the K\"ahler structure
$(\omega,I)$ on $M$. Then it is easy to check that the assumptions
in Theorem \ref{gcarrell} are all satisfied. Thus
$H^k_{\bar{\partial}_2}(M)=0\,\,\,\text{if }\vert k\vert
>\text{dim}Z$. Let $U^k$ be the $-k\sqrt{-1}$ eigenbundle of the generalized complex
structure $\J_2$. It then follows from Example \ref{complexgc} that
$U^i=\bigoplus_{q-p=i} \left(\wedge^{p,q}T^*_{\C}M\right)$ and that
the $\bar{\partial}$ operator associated to the generalized complex
structure $\J_2$ coincides with the usual $\bar{\partial}$
associated to the complex structure $I$. This finishes the proof.

\end{proof}

\section{Calculation of generalized Hodge numbers}

The generalized Hodge theory for compact generalized K\"ahler
manifolds has been established by Gualtieri \cite{Gua04}. Let us
recall some salient points about this theory and refer to
\cite{Gua04} for details.

Let $(M, \mathcal{J}_1, \mathcal{J}_2)$ be a compact generalized
K\"ahler manifold of dimension $2n$. Then
$-\mathcal{J}_1\mathcal{J}_2$, regarded as a positive definite
metric on $TM\oplus T^*M$, induces a Hermitian inner product, the
\textbf{Born-Infeld} inner product, on the space of differential
forms.

Let $\Gamma(U^k)$ be the $-\sqrt{-1}$ eigenspace of the Clifford
action of $\mathcal{J}_1$ on the space of differential forms. The
commuting endomorphism $\mathcal{J}_2$ further decomposes $U^k$ as
\[ \Gamma(U^k)=\Gamma(U^{k,\vert k\vert -n}) \oplus \Gamma(U^{k, \vert k\vert-n+2})\oplus
\cdots \oplus \Gamma(U^{k, n-\vert k\vert}).\]

Thus there is a $(p,q)$ decomposition of the differential forms.
Furthermore, this decomposition is orthogonal with respect to the
\textbf{Born-Infeld} metric, and gives rise to the following
splitting of the exterior derivative:
\[ d=\delta_{+} +\delta_-+\overline{\delta}_++\overline{\delta}_-.\]
Here the differential operators act as follows:

\begin{diagram}
 \Gamma(U^{p-1,q+1}) &        &              &&\Gamma(U^{p+1,q+1})& \\
 &\luTo^{\delta_-}       &   \uDotsto_{\bar{\partial}_2}  &              \ruTo^{\overline{\delta}_+}  &\\
  &    \lDotsto^{\partial_1}   & U^{p,q} &   \,\,  \rDotsto^{\bar{\partial}_1}    &\\
  &\ldTo_{\delta_+}       &\dDotsto^{\partial_2}    &\,\rdTo_{\overline{\delta}_-}&\\
  \Gamma(U^{p-1,q-1}) &             &          &&\Gamma(U^{p+1,q-1})\\
\end{diagram}

Let $d^*$,$\bar{\partial}^*$,$\partial^*$, $\delta^*_{\pm}$ and
$\overline{\delta}^*_{\pm}$ be the adjoint operators of
$d$,$\bar{\partial}$,$\partial$, $\delta_{\pm}$ and
$\overline{\delta}_{\pm}$ with respect to the \textbf{Born-Infeld}
inner product respectively. Then
$d+d^*,\bar{\partial}_{1/2}+\bar{\partial}_{1/2}^*,\partial_{1/2}+\partial_{1/2},\delta_{\pm}+\delta_{\pm}^*,
\overline{\delta}_{\pm}+\overline{\delta}_{\pm}^*$ are all elliptic
operators; moreover,  we have
\[\bigtriangleup_d=2\bigtriangleup_{\partial_{1/2}}=2\bigtriangleup_{\bar{\partial}_{1/2}}=
4\bigtriangleup_{\delta_{\pm}}=4\bigtriangleup_{\overline{\delta}_{\pm}},
\] where  $\bigtriangleup_d$,
$\bigtriangleup_{\partial_{1/2}}$,
$\bigtriangleup_{\bar{\partial}_{1/2}}$,
$\bigtriangleup_{\delta_{\pm}}$,
$\bigtriangleup_{\overline{\delta}_{\pm}}$ are the Laplacians of
$d$, $\partial_{1/2}$, $\bar{\partial}_{1/2}$,
$\bar{\partial}_{1/2}$, $\delta_{\pm}$ and $\bar{\partial}_{\pm}$
respectively. It then follows from the standard Hodge theory for
elliptic operators

\begin{theorem}(\cite{Gua04})\label{gkdecomp}  The cohomology of a compact $2n$
dimensional generalized K\"ahler manifold carries a Hodge
decomposition
\[ H^*(M,\C)=\bigoplus_{\begin{subarray}{1}  \vert p+q \vert \leq n \\
 p+q \equiv 2 (\text{mod} 2)
\end{subarray}}\mathscr{H}^{p,q}, \]
where $\mathscr{H}^{p,q}$ are $\bigtriangleup_d$ harmonic forms in
$\Gamma(U^{p,q})$.

\end{theorem}

 The \textbf{generalized Hodge number} $h^{p,q}$ of a generalized K\"ahler manifold
 $(M,\mathcal{J}_1,\mathcal{J}_2)$
 is then defined to be the complex dimension of $\mathscr{H}^{p,q}$. The case of interest
 is when $(\J_1,\J_2)$ is not the $B$-transform of a genuine K\"ahler structure.
  (c.f. Example \ref{Kahlerexamples}.) In this paper, we refer to such generalized K\"ahler
  structures as \textbf{non-trivial}.  Note that \cite{LT05}
 proposes a general method of constructing non-trivial explicit examples
of generalized K\"ahler structures as the generalized K\"ahler
quotient of the vector space $\C^n$. In the rest of
 this section, we are going to compute the generalized Hodge number for two of these examples.
But let us first recall how to construct non-trivial examples of
generalized K\"ahler manifolds as generalized K\"ahler quotient.

Let $\J$  be a generalized complex structure on a vector space
$V=\C^n$. Let $L \subset V_\C \oplus V_\C ^*$ be the $\sqrt{-1}$
eigenspace of $\J$. Since $L$ is maximal isotropic and $L \cap
\overline{L} = \{0\}$, we can (and will) use the  metric to identify
$L^*$ with $\overline{L}$.

Given $\epsilon \in \wedge^2 L^*$, define $L_\epsilon = \{Y +
\iota_Y\epsilon \mid Y \in L \}.$ Then $L_\epsilon$ is  maximal
isotropic, and $L_\epsilon \cap \overline{L_\epsilon} = \{0\}$ if
and only if the endomorphism
\begin{equation} \label{real index zero condition}  A_{\epsilon}=\left(
\begin{matrix} 1 &  \bar{\epsilon} \\
\epsilon & 1 \end{matrix} \right)  \colon L\oplus \overline{L}
\rightarrow L \oplus \overline{L} \end{equation} is invertible. If
it is invertible, there exists a unique generalized complex
structure $\J_\epsilon$ on $V$ whose $\sqrt{-1}$ eigenspace is
$L_\epsilon$. Note that $A_\epsilon$ is always invertible for
$\epsilon$ sufficiently small.

Now let $(\J_{\omega}, \J_I)$ be the generalized K\"{a}hler
structure on $V=\C^n$ which is induced by the standard genuine
K\"ahler structure $(\omega,I)$ . Let $L_1$ and $L_2$ denote the
$\sqrt{-1}$ eigenspaces of $\J_{\omega}$ and $\J_I$ respectively.
Then $L_1 =\left(L_1 \cap L_2 \right) \oplus \left( L_1 \cap
\overline{L_2} \right)$ and $L_2 =\left(L_1 \cap L_2 \right) \oplus
\left( \overline{L_1}\cap L_2 \right)$. Thus $\epsilon \in
C^\infty(\wedge^2 \overline{L_2})$ fixes $\J_{\omega}$ if and only
if $\epsilon$ takes $L_1 \cap L_2$ to $L_1 \cap \overline{L_2}$,
i.e., if and only if $\epsilon$ is an element of
$C^\infty\left((\overline{L_1} \cap\overline{L_2}) \otimes (L_1 \cap
\overline{L_2})\right)$.

Henceforth we will assume that $\epsilon \in
C^\infty\left((\overline{L_1} \cap\overline{L_2}) \otimes (L_1 \cap
\overline{L_2})\right)$. So it will fix $\J_{\omega}$ and deform
$\J_I$ to a new generalized almost complex structure $\J_{\epsilon}$
on a bounded region of $\C^n$.  The following lemma gives a simple
condition which guarantees that $L_\epsilon$, the $\sqrt{-1}$
eigenbundle of $\J_{\epsilon}$, is closed under the Courant bracket,
and hence that $(\mathcal{J}_{\omega},\J_\epsilon)$ is a generalized
K\"ahler structure.

\begin{lemma} \label{closedness2}
Assume that there exists a subset $I \subset (1,\ldots,n)$ so that
$$\epsilon  =
\sum_{i,j \in I}  F_{ij} (z) \frac{\partial}{\partial z_i} \wedge
\frac{\partial}{\partial z_j} + \sum_{i,j  \in I} F_{ij}(z)
d\overline{z_i} \wedge d\overline{z_j}.$$ If $F_{ij}$ is holomorphic
and $\dfrac{\partial{F_{ij}}}{\partial z_k} = 0$ for all $i, j$ and
$k \in I$, then $L_\epsilon$ is closed under the Courant bracket.
\end{lemma}

Consider the action of a compact connected torus $T$ on $V=\C^n$
which preserves the generalized K\"ahler structure
$(\J_{\omega},\J_{\epsilon})$. A \textbf{generalized moment map} for
$(\mathcal{J}_{\omega}, \mathcal{J}_{\epsilon})$ is the generalized
moment map for the generalized complex structure
$\mathcal{J}_{\omega}$ which coincides with the usual moment map for
the symplectic structure $\omega$ in this context. Let $\t$ be the
Lie algebra of the torus $T$ and $f: M \rightarrow \t^*$ be the
generalized moment map for the torus action. If $a \in \t^*$ so that
$T$ acts freely on $f^{-1}(a)$, then $M_a=f^{-1}(a)\diagup S^1$ is
defined to be the \textbf{generalized K\"ahler quotient} of the
torus action at the level $a$. There is a naturally defined
generalized K\"ahler structure
$(\widetilde{\mathcal{J}}_1,\widetilde{\mathcal{J}}_2)$ on $M_a$;
moreover, for any $m \in f^{-1}(a)$ we have
\begin{equation} \label{typeformula}\begin{split}
& \text{type}(\widetilde{\J}_{\omega})_{[m]}=\text{type}(\J_{\omega})_m=0,\\
&\text{type}(\widetilde{\J}_{\epsilon})_{[m]}=\text{type}(\J_{\epsilon})_m-\text{dim}(T)+2\text{dim}(\t_M\cap
\pi(L_{\epsilon}))_m,\end{split}\end{equation} where $\pi:
T_{\C}\C^n\oplus T_{\C}^*\C^n \rightarrow T_{\C}\C^n$ is the
projection map, $\t_M$ is the distribution of fundamental vector
fields generated by the torus action, and $L_{\epsilon}$ is the
$\sqrt{-1}$ eigenbundle of the generalized complex structure
$\mathcal{J}_{\epsilon}$. Since a $B$-transform always preserves the
type of a generalized complex structure,
$(\widetilde{\J}_{\omega},\widetilde{\J}_{\epsilon})$ is a
non-trivial generalized K\"ahler structure if and only if
$\text{type}(\widetilde{\J}_{\epsilon})_x \neq n-\text{dim}(T)$ at
some point $x$ of $M_a$.

\begin{example} \label{ghodge1}({\bf  $\CP^n$ for $n \geq 3$})

 We now construct a non-trivial
 generalized K\"ahler structure on $\CP^n$ for $n\geq 3$
and compute its generalized Hodge number.

Let $S^1$ act on $C^{n+1}$ via
\[ \lambda (z_0, \cdots,z_n)=(\lambda z_0,\cdots,\lambda z_n).\]
Note that this action preserves the K\"ahler structure $(\omega,I)$
with a moment map given by
\[\varPhi(z)=\sum_i\dfrac{1}{2}\vert z_i\vert^2;
\] moreover, the $S^1$ action on $\varPhi^{-1}(1)$ is free and the
reduced space is exactly $\CP^n$.

Let \[ \epsilon = z_0z_1\dfrac{\partial}{\partial z_2}\wedge
\dfrac{\partial}{\partial z_3} +z_0z_1d\overline{z}_2\wedge
d\overline{z}_3.\] If necessary, multiplying $\epsilon$ by a
sufficiently small positive constant, then $\epsilon$ will deform
$(\mathcal{J}_{\omega}, \mathcal{J}_I)$ to a new generalized almost
K\"ahler structure on the bounded open set $O=\{\varPhi^{-1}(z)
<2\}$ so that $\text{type}(\J_{\epsilon})_z$ is $n+1$ if $z_0z_1=0$
 and is $n-1$ otherwise. Since $z_0z_1$ is holomorphic and
$\dfrac{\partial z_0z_1}{\partial z_2}=\dfrac{\partial z_0z_1
}{\partial z_3}=0$, by Lemma \ref{closedness2}
$(\mathcal{J}_{\omega},\mathcal{J}_{\epsilon})$ is in fact a
generalized K\"ahler structure on $O$.

Since $\epsilon$ is $S^1$ invariant,
$(\mathcal{J}_{\omega},\mathcal{J}_{\epsilon})$ is also $S^1$
invariant with a generalized moment map $\varPhi$. So there is a
naturally defined generalized K\"ahler structure
$(\widetilde{\J}_{\omega},\widetilde{\J}_{\epsilon})$  on the
quotient space $\CP^n=\varPhi^{-1}(1)\diagup S^1$. Note that  the
fundamental vector generated by the above $S^1$ action on $\C^{n+1}$
is
$$X =\dfrac{\newi}{2}\sum_i \left( z_i \dfrac{\partial }{\partial z_i}
-\overline{z}_i \dfrac{\partial }{\partial \overline{z}_i} \right)
\,,$$ and hence $X$ does not lie in $\pi(L_\epsilon)$ at any point
of $\C^{n+1}$, where $L_\epsilon$ is the $\sqrt{-1}$ eigenbundle of
$\mathcal{J}_\epsilon$. It follows from (\ref{typeformula})
$\text{type}(\widetilde{\J}_{\omega})_{[z]}=0$ for all $[z] \in
\CP^n$, whereas $\text{type} (\widetilde{\J}_{\epsilon})_{[z]}=n$ if
$z_0z_1=0$ , otherwise $\text{type}
(\widetilde{\J}_{\epsilon})_{[z]}=n-2$. So by the preceding
discussion the generalized complex structure
$(\widetilde{\J}_{\omega},\widetilde{\J}_{\epsilon})$ on $\CP^n$ is
non-trivial.

Next consider the $S^1$ action on $\C^{n+1}$ defined by
\[\lambda (z_0,z_1,z_2, z_3,\cdots ,z_i,\cdots, z_n)=
(\lambda z_0,\lambda^5 z_1, \lambda^2z_2, \lambda^4 z_3,\cdots
,\lambda^{2i} z_i ,\cdots ,\lambda^{2n}z_n).\] It is easily seen
that this action preserves the standard K\"ahler structure and the
deformation $\epsilon$; furthermore it commutes with the standard
$S^1$ action on $\C^{n+1}$ and so descends to an action on $\CP^n$
which preserves the quotient generalized K\"ahler structure
$(\widetilde{\J}_{\omega},\widetilde{\J}_{\epsilon})$. It is
equivariantly formal since the action of any compact Lie group on
$\CP^n$ is. Moreover, it is easy to check that this action has $n+1$
isolated fixed points:
\[ [(1,0,\cdots,0)],[(0,1,\cdots,0)],\cdots,[(0,0,\cdots,1)].\]
Observe at the tangent space of each fixed point $x$ the generalized
K\"ahler structure $\widetilde{\J}_{\epsilon,x}$ is induced by a
complex structure. It then follows from Theorem \ref{gcarrell} that
\[H^i_{\bar{\partial}_{\widetilde{\J}_{\epsilon}}}(\CP^n)=0,\,\,\,\text{if}\,\,i\neq
0.\] Thus the generalized Hodge number $h^{p,q}=0$ if $q\neq 0$. In
addition, we have \[\begin{split}
h^{p,0}&=\text{dim}H^p_{\bar{\partial}_{\widetilde{\J}_{\omega}}}(\CP^n)
\,\,\,\,(\text{By Theorem \ref{gkdecomp}}.)
\\&=\text{dim}H^{n+p}(\CP^n)\,\,\,(\text{By
Example \ref{symplecticgc}}.)\\&=
\begin{cases} 1 \,\,\,\,\,\text{if $n+p$ is even;}\\
0\,\,\,\,\text{if $n+p$ is odd.}\end{cases}
\end{split}\]

\end{example}

\begin{example}{\bf $\CP^n$ blown up at one point for $n\geq 3$}

In this example we  construct a non-trivial generalized K\"ahler
structure on the blow up of $\CP^n$ at one point and compute its
generalized Hodge number.

Let a two dimensional torus $T$ with Lie algebra $\t$ act on
$\C^{n+2}$ by
\[(\alpha,\beta)(z_1,\cdots,z_n,z_{n+1},z_{n+2})=(\alpha \beta
z_1,\cdots,\alpha \beta z_n,\alpha z_{n+1},\beta^{-1}z_{n+2})\] with
moment map
\[f(z_1,\cdots,z_n,z_{n+1},z_{n+2})=(\vert z_1\vert^2+\cdots+\vert z_n\vert^2+\vert
z_{n+1}\vert^2, \vert z_1\vert^2+\cdots+\vert z_n\vert^2-\vert
z_{n+2}\vert^2).\]   Then there exists some $\xi \in \t^*$ so that
$M_{\xi}=f^{-1}(\xi)/ T^2$ is equivariantly symplectomorphic to
$\CP^n$ blown up at a fixed point.

Define
\[\epsilon =(z_1z_2z_{n+2}) \dfrac{\partial}{\partial z_n} \wedge
\dfrac{\partial}{\partial z_{n+1}} +(z_1z_2z_{n+2})
d\overline{z}_n\wedge d \overline{z}_{n+1}\,.\] It is clear from
construction that $\epsilon$ is $T$ invariant. As explained in
\cite{LT05}, $\epsilon$ deforms the standard K\"ahler structure
$(\J_{\omega},\J_{I})$ to a $T$ invariant generalized K\"ahler
structure $(\J_{\omega},\J_{\epsilon})$; moreover,
$(\J_{\omega},\J_{\epsilon})$ descends to a non-trivial generalized
K\"ahler structure
$(\widetilde{\J}_{\omega},\widetilde{\J}_{\epsilon})$ on the reduced
space $M_{\xi}$, i.e., $\CP^n$ blown up at a fixed point.

Define $S^1$ action on $\C^{n+2}$ by \[ \alpha
(z_1,\cdots,z_n,z_{n+1},z_{n+2})=(\alpha^{\lambda_1}z_1,\cdots,\alpha^{\lambda_n}z_n,\alpha^{\lambda_{n+1}}z_{n+1},\alpha^{\lambda_{n+2}}z_{n+2}),\]
where $\lambda_1,\cdots,\lambda_n,\lambda_{n+1},\lambda_{n+2}$ are
rational numbers so that
\begin{enumerate}\item [a)] $\lambda_1,\cdots,\lambda_n$ are
distinct from each other;
\item [b)]
$\lambda_1+\lambda_2+\lambda_{n+2}=\lambda_n+\lambda_{n+1}$;
\item [c)] $\lambda_i \neq \lambda_{n+1}-\lambda_{n+2}$,
$i=1,2,\cdots,n$. \end{enumerate}

This $S^1$ action preserves the standard K\"ahler structure
$(\J_{\omega},\J_I)$ on $\C^{n+2}$ and the deformation $\epsilon$;
furthermore it commutes with the $T$ action and therefore descends
to a $S^1$ action on $\CP^n$. It is easy to check that the induced
$S^1$ action on $\CP^n$ is Hamiltonian and so is equivariantly
formal. Moreover, by construction the induced action on $\CP^n$ has
only finitely many fixed points such that the restriction of the
generalized complex structure $\widetilde{\J}_{\epsilon}$ to the
tangent space of each fixed point is a complex structure.

 It then follows from Theorem \ref{gcarrell} that
\[H^i_{\bar{\partial}_{\widetilde{\J}_{\epsilon}}}(M)=0,\,\,\,\text{if}\,\,i\neq
0.\] Thus the generalized Hodge number $h^{p,q}=0$ if $q\neq 0$. In
addition, we have \[\begin{split}
h^{p,0}&=\text{dim}H^p_{\bar{\partial}_{\widetilde{\J}_{\omega}}}(M)\,\,\,\,(\text{By
Theorem \ref{gkdecomp}}.)\\&=\text{dim}H^{n+p}(M)\,\,\,(\text{By
Example \ref{symplecticgc}}.)
\end{split}\]
\end{example}

\end{document}